%%%%%%%%%%%%%%%%%%%%%%%%%%%%%%%%%%%%%%%%%%%%%%%%%%%%%%%%%%%%%%%%%%%%%%%%%%%
%% Jaffe, Arthur; Quinn, Frank
%% 
%% Response to comments on ``Theoretical Mathematics''
%% 
%% publ:  Bull. Amer. Math. Soc. (N.S.) 30(1994) no. 2
%% pp:    208-211
%% type:  Research-Expository Paper    markup: amstex    file size: 13K
%% 
%% copyright: American Math. Society copyright; see end of article
%% 
%% Include files necessary for this article: bull-ppt.tex
%% 
%%%%%%%%%%%%%%%%%%%%%%%%%%%%%%%%%%%%%%%%%%%%%%%%%%%%%%%%%%%%%%%%%%%%%%%%%%%
\input amstex
\documentstyle{amsppt}
\input bull-ppt

\topmatter
\cvol{30}
\cvolyear{1994}
\cmonth{April}
\cyear{1994}
\cvolno{2}
\cpgs{208-211}
\title Response to comments on ``Theoretical Mathematics''
\endtitle

\author Arthur Jaffe and Frank Quinn\endauthor
\endtopmatter

\document
\heading Reactions\endheading
We would like to take this
opportunity to thank the many people who have
sent notes and made comments complimenting our work.
Most reaction we have received
from the mathematics community
has been very favorable.
In this volume the editor has focused on divergent views,
for instance by soliciting replies from people named in the
original article. We would also like to express our 
appreciation
for the care and thought that these authors have
put
into these replies, and
for the opportunity to respond to them in the same spirit.
We received some other reactions from physicists,
see the response of Friedan,
but we
would like to emphasize that the article was not addressed
to physicists doing physics. A style of working which may be
counterproductive in mathematics may be totally appropriate
in physics; this is an entirely different matter.
Our article was addressed to mathematicians and to
physicists who want their work considered as mathematics.
In any case, a dialog has begun.

\heading Terminology \endheading
Perhaps the most frequent objection to our paper is
not to the substance, but rather to our use of the term
``theoretical''.
A primary objective of our paper was to welcome and 
celebrate
speculation.
To do this, we required a term with both dignity and 
resonance.
``Speculator'' did not seem to
be such a word, nor was ``conjecturer''. We felt that
``theoretical mathematician'', with its evocation of 
``theoretical
physicist'' but with its clear identification with 
mathematics,
satisfied these criteria.
We agree there is a discrepancy between our use of
``theoretical'' to mean ``speculative'' and the 
well-established
mathematical
use of ``theory'' to refer to a coherent body of knowledge.
%Our use is closer to the common-language use,
But our main
motivation was that we could not find
another term with similar qualities.

\heading The value of conjecture \endheading
Most people correctly read our paper as positive about
speculation and designed to enhance it. Some, like Mac Lane,
were concerned that we had gone too far.  A few, like
Atiyah, read it as hostile to conjecture. We certainly agree
with him that suppression of speculation would condemn
mathematics to ``an arthritic old age''. But we do not see
that calling it by its true name, or explicitly 
acknowledging
that it is incomplete, would have this effect.

Along the same lines, we appreciated Thom's remark that
``rigor'' reminds one of the phrase ``rigor mortis''.  There
is a very real sense in which mathematics which has been
successfully rigorized is dead, and the real life is, as
Atiyah suggests, in speculation.  We might think of
mathematics as a tree: only the leaves and a thin layer
around the trunk and branches are actually ``alive''. But 
the
tree is supported by a large literature of published 
``dead''
wood. Our position is that material should be really sound
before being allowed to die and to be incorporated into 
the wood. In
this image the ``self-correcting nature of science'' is that
if too much rotten wood is published then the branch falls
off and the tree starts over.

\heading Who cares?\endheading
It is evident that few of the responders have themselves 
been
bothered by the problems we describe.  Where, then, is the 
problem?
Gray suggests ``all those who work away from the main 
centers of
research are disadvantaged.''
``Away'' here should include  time as well as space,
to encompass those who use the mathematical works
from previous generations.
Further,
there are mathematicians of
such power and insight that they
have no difficulty evaluating new work, extracting value, 
and
discarding errors and chaff.
Everyone else is disadvantaged.
The large majority of mathematicians, as well as 
nonmathematicians
who use mathematics, rely on the literature.
This issue may be one of the
rare occasions when the leadership offered by
some of the mathematical elite is inappropriate
for average mathematicians.

\heading The effect on students\endheading
We wrote that students might be adversely affected by the
problems we discussed.  Atiyah and Uhlenbeck, among others,
seemed to interpret this as suggesting that students should
be warned away from areas of speculation (and excitement).
This was not our intent.  We both sought out such areas as
students and have no regrets.  Rather, it is a question of
generally accepted standards, goals, and role models.
We suggest it is a dangerous
thing for a student who does not understand the special 
character of
speculative work to set out to emulate it. Unless the 
student has
extraordinary talent, this will work out poorly.

\heading Research Announcements\endheading
Many people accept our premises, arguments, and other
conclusions but balk at the conclusion that published
research announcements prior to the existence of an
article are unwholesome. We reiterate that we
make no objection to the circulation of research
announcements to inform or to
establish priority.  We object only
to their official publication. This conclusion came as a
shock to us too: at the time one of us was Research
Announcements editor for the {\it Bulletin\/} and had a 
strong
commitment to the format. But we found that the 
understanding we
developed compelled this conclusion.

The main argument in support of announcements is that they
have proved themselves in practice. But the world is 
changing
rapidly, and the experiences of the past must be 
extrapolated
into the future very carefully. In particular electronic
communications are certain to change the nature of
publication. We see announcements as being far more
problematic in this new world than they were in the old.

\heading Specificity and the naming of names\endheading
We were taken to task by several respondents for limiting 
our
discussion of the mathematics-science interaction to 
theoretical
physics. Certainly it is true that there are
profound interactions in other areas, but we chose to focus
our discussion for three reasons.  First, the activity in 
the
specific area is unusual in quantity and quality.  As
Atiyah puts it, ``It involves front-line ideas both in
theoretical physics and in geometry.''  Second, we both have
direct experience with interactions in this area.
Third, we felt it important to be specific. An
argument with some detail in a specific context seems to us
to be better than vague generalities or widely scattered
examples.
We believe that the conclusions extracted from this 
particular case
apply more generally.

We apologize to those
who experienced discomfort because we used the names of 
living
mathematicians.
Unfortunately some people took our mentioning of names as a
personal attack.  It was certainly not intended this way.
We only mentioned
exceptionally well-known, highly regarded mathematicians 
with the
respect of the entire community.  We thought that these 
people,
whose contributions and reputations are beyond question,
might be
discussed at a level which transcended personalities.

\heading Reply to Thurston\endheading
Thurston wrote an extremely thoughtful response to our 
paper,
and we hope the readers appreciate it as much as we do.
He carefully and eloquently articulated a vision of
mathematics, in some ways quite different from ours. But 
as he points
out, the way in which a question is phrased can greatly 
influence the
answer.  We feel that his carefully crafted questions do 
lead to
important insights, but they also channel discussion away 
from
our particular concerns.

His lead question is, How do mathematicians advance human
understanding of mathematics? This is followed by, How do 
people
understand mathematics? These introduce an analysis of
the limitations and needs of real people and make a 
convincing case that
the formal language of mathematical papers is poorly 
adapted to
these needs. He argues that understanding would be better 
achieved
through other avenues, including a relaxation of the 
emphasis on
rigor.

But understanding comes on many different levels. We believe
that Thurston describes the level of understanding a 
teacher might want to
instill in students. It is sufficient for appreciation but 
usually not for
active application or further development. Traditional 
papers are,
as he says, a poor way to introduce a subject, but they 
are the
whetstones which give the final edge to mastery. And they 
are at
least a last resort in the learning process. Thurston
himself may obtain satisfactory understanding through 
informal
channels, but he
is a mathematician of extraordinary power and  should be 
very
careful about extrapolating from his experiences to the 
needs of
others.

Rigorous papers serve other functions besides education, 
including
the way we determine what is true and reliable. And proofs 
provide a
source of techniques for other problems and clues to new
and unsuspected phenomena. Thurston relates his 
introduction to the
social aspects of knowledge and says, ``This knowledge was 
supported
by written documents, but the written documents were not 
really
primary. \dots Andrew Wiles's proof of Fermat's Last 
Theorem serves as a
good illustration of this,\dots The experts quickly came 
to believe
that his proof was basically correct on the basis of 
high-level ideas,
long before details could be checked.'' But in their 
pronouncements
the experts were careful to note that only the details 
could be
authoritative. And sure enough, at this time (December 1993)
the details are still inconclusive.
The result is still uncertain, and
we might have more to learn from this wonderfully
fruitful quest. But in any case the written documents 
really are
primary.

Later Thurston asserts ``\!\dots our strong communal 
emphasis on
theorem-credits has a negative effect on mathematical 
progress,''
since it obscures and impedes the team aspect of 
mathematics. This
is one of our  points too; theoretical mathematicians in
particular must recognize they are part of a team, not the 
whole
show, and must be willing to share credit.

Clearly Thurston has identified real weaknesses and
needs in the educational, social, and communication 
aspects of
mathematics. But his analysis is not the whole story. He 
has steered
the discussion away from what seem to us to be key issues, 
particularly
the importance of ``truth in advertising'' in the 
publication
of research and the ultimate goal of proof in the
written record. Curiously, we feel that his analysis may 
apply better
to sciences other than mathematics, where the tools for
establishing reliability
are far less effective than in mathematics.  But it seems
to us that even in mathematics when his conclusions are 
limited to
education rather than research, they are complementary to 
and compatible
with ours.

\enddocument